\documentclass[11pt]{amsart}
\usepackage{amssymb}
\usepackage{verbatim}
\input xypic
\xyoption{all}

\newcommand{\BZ}{{\mathbb{Z}}}

\newcommand{\BC}{{\mathbb{C}}}
\newcommand{\BF}{{\mathbb{F}}}

\newcommand{\BP}{{\mathbb{P}}}

\newcommand{\gb}{\beta}

\newcommand{\gc}{\gamma}

\newcommand{\gs}{\sigma}
\newcommand{\gS}{\Sigma}
\newcommand{\gO}{\Omega}

\newcommand{\ga}{\alpha}
\newcommand{\gt}{\tau}

\newcommand{\gth}{\theta}

\newcommand{\gT}{\Theta}

\newcommand{\fB}{{\mathfrak{B}}}
\newcommand{\fC}{{\mathfrak{C}}}

\newcommand{\fR}{{\mathfrak{R}}}

\newcommand{\cM}{{\mathcal{M}}}

\newcommand{\cO}{{\mathcal{O}}}

\newcommand{\ti}[1]{\tilde{#1}}

\newcommand{\ol}[1]{\overline{#1}}

\newcommand{\Pic}{\mathrm{Pic}}
\newcommand{\jac}{\mathrm{Jac}}

\newcommand{\prym}{\mathrm{Prym}}
\newcommand{\gal}{\mathrm{Gal}}
\newcommand{\Sym}{\mathrm{Sym}}

\newcommand{\norm}{\mathrm{Norm}}

\newcommand{\Id}{\mathrm{Id}}

\newcommand{\Cy}[1]{\BZ/#1\BZ}
\newcommand{\hra}{\hookrightarrow}
\newcommand{\ra}{\longrightarrow}

\newcommand{\sm}{\smallsetminus}
\swapnumbers
\theoremstyle{plain}
\newtheorem{lma}{Lemma}[section]
\newtheorem{thm}[lma]{Theorem}
\newtheorem{prp}[lma]{Proposition}
\newtheorem{qst}[lma]{Question}
\newtheorem{cor}[lma]{Corollary}

\theoremstyle{definition}
\newtheorem{dfn}[lma]{Definition}
\newtheorem{rmr}[lma]{Remark}
\newtheorem{ntt}[lma]{Notation}

\newtheorem{dsc}[lma]{}

\begin{document}

\title[reconstruction of quartics]
{Any smooth plane quartic can be reconstructed from its bitangents}
\thanks{The author was partially supported by Israel-US BSF grant 1998265.}
\author{David Lehavi}
\address{Department of Mathematics\\ 
The Ohio State University \\
100 Math Tower \\
231 West 18th Avenue\\ 
Columbus, OH 43210-1174, USA}
\email{dlehavi@math.ohio-state.edu}
\date{\today}
\begin{abstract}
In this paper, we present two related results on curves of genus 3.
The first gives a bijection between the classes of the following objects:
\begin{itemize}
\item Smooth non-hyperelliptic curves $C$ of genus 3, with a choice of an
element $\ga\in\jac(C)[2]\sm\{0\}$, such that the cover 
$C\ra|K_C+\ga|^*$
does not have an intermediate factor; up to isomorphism.
\item Plane curves $E,Q\hra\BP^2$ and an element in
$\gb'\in\Pic(E)[2]\sm\{0\}$, where $E,Q$ are of
degrees 3,2, the curve $E$ is
smooth and $Q,E$ intersect transversally ; up to
projective transformations.
\end{itemize}
We discuss the degenerations of this bijection, and give an interpretation of
the bijection in terms of Abelian
varieties. Next, we give an application of this correspondence:
An {\em explicit} proof of the reconstructability of {\em any} smooth plane
quartic from its bitangents.
\end{abstract}
\maketitle

%
\section{Introduction}\label{Sintro}
%
\begin{dsc}
It is a classical result that {\em any smooth plane quartic, has exactly 28 bitangents.}
The properties of these bitangents were a popular research subject from the
second half of the 19th century, until the 1930's (Cayley, Hesse, Pl\"ucker,
Schottky, Steiner, Weber). A natural question to ask is:
\end{dsc}
\begin{qst}
Can a smooth plane quartic be
reconstructed from data related to its 28 bitangents ?
\end{qst}
\begin{dsc}
Various partial answers were given. Aronhold proved (see \cite{KW} p. 783 for
a proof ) that smooth plane quartics can be reconstructed from certain
7-tuples of bitangents called ``Aronhold system''.
Recently Caporaso \& Sernesi proved in \cite{CS1} that a {\em generic}
curve of degree $d\geq 4$ is defined by its bitangents, and in
\cite{CS2} that a {\em generic} curve of genus $g>3$ is
determined by the image of its odd theta-characteristics under the
canonical embedding. The proof in both papers is identical in the
plane quartic case. It involves degenerations. Thus to the best of this
authors knowledge, can not be made either non-generic or constructive.

The main theorem of this paper is:
\end{dsc}
\begin{thm}\label{Tbitan}
Any smooth plane quartic $C$ can be reconstructed from its bitangents.
\end{thm}
\begin{dsc}
The proof is explicit - one can derive an explicit formula from the proof.
\end{dsc}
\begin{ntt}\label{NCa}
The following notation is kept throughout the paper:
We fix the ground field as $\BC$ (see Remark \ref{Rfields} regarding other fields).
We consider pairs $(C,\ga)$, where $C$ is a smooth curve of genus 3,
and
\[
  \ga\in\jac(C)[2]\sm\{0\},
\]
such that the linear system $|K_C+\ga|$ does not have base points.
\end{ntt}
\begin{dsc}
In the course of the proof (in Theorem \ref{Tpc}) we identify three moduli
spaces:
\end{dsc}
\begin{thm}\label{Tpc}
There is a coarse moduli space isomorphism between the moduli of the following
sets of data:
\begin{itemize}
  \item Pairs $(C,\ga)$ (see Notation \ref{NCa}) such that the cover 
$C\ra|K_C+\ga|^*$
does not have an intermediate factor; up to isomorphisms.
  \item Pairs $(Y/E,\gb')$ such that $Y/E$ is a ramified double cover 
    of irreducible smooth curves of genera $4,1$ respectively,
   $\gb'\in\Pic(E)[2]\sm\{0\}$; up to isomorphisms.
  \item Plane curves configurations $E,Q\hra\BP^2$, of degrees 3,2
   respectively and an element $\gb'\in\Pic(E)[2]\sm\{0\}$, such
   that the curve 
   $E$ is smooth, and the curves $Q,E$ intersect transversally;
   up to projective transformations.
\end{itemize}
\end{thm}
\begin{dsc}
The generic case of Theorem \ref{Tpc} is classical (see \cite{Co}
section 48, \cite{BC}, \cite{Dol}, \cite{Dol2}).
The view presented here has the following advantages over the classical view:
\begin{itemize}
\item We describe some degeneration loci:
The hyperelliptic locus is described in
Theorem \ref{TChee} and Proposition \ref{Psing}. The locus for which the
cover $C\ra|K_C+\ga|^*$ has an intermediate factor is described in Section
\ref{Ssing}.
\item We use the fact that Theorem \ref{Tpc} is {\em not} generic, and
the description of the degeneration loci, to prove Theorem \ref{Tbitan}.
\item In the last part of section \ref{Sprym} we give an Abelian varieties
interpretation of Theorem \ref{Tpc}, relating the Jacobians of the curves
$C,Y,E$. This interpretation is useful in evaluating
Abelian integrals (see \cite{Le}).
\end{itemize}
In the proof of Theorem \ref{Tbitan} we use the classification of finite
2-transitive groups. See Remark \ref{Rtrans} for details.
\end{dsc}
\subsection*{Acknowledgments}
I thank Ron Donagi, Hershel Farkas, Assaf Libman, Shahar Mozes, 
Christian Pauly and Christophe Ritzenthaler for the
illuminating conversations we had during my work on this paper.
Igor Dolgachev helped me put this work in the correct historical
perspective. Dan Abramovich and Tamar Ziegler helped me bring it to a readable
form. Most of all I thank Ron Livn\'e for introducing me to genus 3 problems,
for his insights, and for his enormous patience with my arguments.
%
\section{Technical Background}
%
\begin{dfn}\label{DPmum}
Let $X$ be a curve of genus $2g-1$, $i:X\ra X$ a fixed point free involution.
Define $Y:=X/i$. Note that the genus of $Y$ is $g$. Consider the morphism
$\norm:\jac(X)\ra\jac(Y)$.
Mumford showed (see \cite{Mum3} p. 331) that
\[
  \ker(\norm(i))=P_0\cup P_1,
\]
where the variety $P_0\subset\jac(X)$ is a principally polarized Abelian
variety, of dimension $g-1$, and the variety $P_1$ is a shift of $P_0$.
Moreover:
\[
  \gT_X\cap P_0=2\Xi,
\]
where $\Xi$ is a principal polarization on $P_0$.
The variety $P_0$ is called the {\em Prym variety} of the double cover $X\ra Y$.
\end{dfn}
\begin{dfn}\label{DPbea}
Beauville generalized definition \ref{DPmum} (see \cite{Be}) to the class of
singular curves $X$ such that:
\begin{itemize}
  \item $X$ has only ordinary double points.
  \item The only fixed points of $i$ are the singular points of $X$.
\end{itemize}
These covers are called {\em admissible covers}.
\end{dfn}
\begin{dfn}
The Prym variety is defined for any double cover, as the 0 component of
$\ker(\norm)$. However, in this case, the Prym variety is not necessarily
principally polarized.
\end{dfn}
\begin{lma}[\cite{DL}, Lemma 1]\label{LDLisogeny}
Let $\pi:\ti{C}\ra C$ be an admissible double cover, and let
 $\nu\pi:\nu\ti{C}\ra \nu C$ be
its partial normalization at $r>1$ points $x_1\ldots x_r\in C$. Let g be the
arithmetic genus of the partial normalization $\nu C$, so the arithmetic
genus of $C$ is $g+r$. Then $\prym(\ti{C}/C)$ has a principal polarization,
$\prym(\nu\ti{C}/\nu C)$ has a polarization of type $2^g 1^{r-1}$, and the
pullback map
\[
  \nu^*\prym(\ti{C}/C)\ra\prym(\nu\ti{C}/\nu C)
\]
is an isogeny of degree $2^{r-1}$.
\end{lma}
\begin{ntt}\label{Nsing}
Given a ramified double cover $\pi:\ti{C}\ra C$, with ramification points 
$x_1\ldots x_r$, choose a partition to pairs on a subset of $x_1\ldots x_r$.
Denote these pairs by $x_{1,1},x_{1,2},\ldots,x_{j,1},x_{j,2}$. We define
the {\em singularization} on these pairs to be the double cover $\ti{C'}\ra C'$
where:
\[
  \begin{aligned}
   C'&:=C/\{x_{1,1}\sim x_{1,2},\ldots, x_{j,1}\sim x_{j,2}\}  \\
   \ti{C'}&:=\ti{C}/\{\pi^{-1}(x_{1,1})\sim \pi^{-1}(x_{1,2}),\ldots, 
                      \pi^{-1}(x_{j,1})\sim \pi^{-1}(x_{j,2})\}.
  \end{aligned}
\]
\end{ntt}
\begin{dsc}\label{DMum}
We will use the symplectic structure on $\jac(C)[2]$ (for a construction of
this structure, in an Abelian varieties context, see \cite{Mum1} p. 183).
Let $C$ be a curve, and let $L$ be a theta characteristic of $C$. Then
\[
  \begin{aligned}
  q:\jac(C)[2]&\ra     \BF_2 \\
           \ga&\mapsto h^0(L+\ga)+h^0(L)\mod 2
  \end{aligned}
\]
is a quadratic form over $\BF_2$. The Weil pairing on $\jac(C)[2]$ is given by
\[
  \langle\ga,\gb\rangle=q(\ga+\gb)-q(\ga)-q(\gb).
\]
Mumford showed (see \cite{Mum2} p. 184) that the form $q$ is locally constant on
families.
\end{dsc}
\begin{ntt}
Denote by $K(V)$ the total field of functions of a scheme $V$.
\end{ntt}
\begin{lma}[The monodromy argument]\label{Lmonodromy}
Let $X$ be an irreducible scheme, $Z$ a finite irreducible Galois extension.
Let $U,V$ be intermediate extensions.
The following decomposition (to irreducible components) holds:
\[
  K(V\times_X U)=\bigoplus_{U,U' \text{ are Galois conjugates}}K(V)\vee K(U').
\]
\end{lma}
\begin{proof}
This follows from the Galois correspondence.
\end{proof}
\begin{dsc}
Typically, we will use the monodromy argument (Lemma \ref{Lmonodromy}) in the
following setting:
\begin{itemize}
  \item $X$ parameterizes curves of genus 3 with some open constraint.
  \item $Z$ adds full level 2 structure to $X$.
\end{itemize}
The monodromy group in this case is the group $SP_6(2)$ (see \cite{Har} 
Result (ii), p. 687).
\end{dsc}
\begin{rmr}\label{Rfields}
Most of the techniques we use are algebraic. We list here the few parts in
which we use analytic arguments, or arguments that use high enough
characteristic of the ground field $\BF$:
\begin{itemize}
  \item Donagi's dictionary for the trigonal and bigonal constructions (see
    \cite{Do} p. 74 and p. 68-69 respectively) requires
    $\mathrm{char}(\BF)\neq 2,3$
  \item In the proof of Theorem \ref{Tbitan}, we use
    Harris's theorem on the Galois group acting on the bitangents of a plane
    quartic (see \cite{Har}, Result (ii), p. 687. Note that the group
    $SP_6(2)$ is denoted there as $O_6(2)$). Harris's proof is analytic.
    In \cite{CCNPW} p. 46, there is a sketch of an algebraic calculation
    of this Galois group over
    a hyperelliptic curve (with odd theta characteristics, instead of
    bitangents, and for high enough or zero characteristic). For high
    enough characteristic, one can try to consider the hyperelliptic case
    as a degeneration of the general case, to get an alternative calculation
    of the Galois group.
\end{itemize}
\end{rmr}
\begin{rmr}\label{Rtrans}
In the proof of Theorem \ref{Tbitan}, we use the classification
of 2-transitive groups, which uses the classification of simple groups.
This is used only to show (in Proposition \ref{PsubS28}) that the group
$SP_6(2)$ is maximal in the group $A_{28}$ under the embedding
$SP_6(2)\hra S_{28}$.
One can try to give a computerized proof of this maximality property 
by adding elements to $SP_6(2)\hra S_{28}$. Construction of the embedding is
given in either \cite{DM} p. 246, \cite{Har} p. 705 or \cite{CCNPW} p. 46.
For generators of this embedding, see \\
\verb|http://www.mat.bham.ac.uk/atlas/v2.0/clas/S62/|
\end{rmr}
%
\section{representing level-2 information with a double cover}\label{SCZ}
%
\begin{dsc}\label{trigo}
Our objective in this section is to define curves $W,Z,X,E,F,Y$,
morphisms $f,g,p_X,p_F,p_Y,q_Y,q_F,q_X$ and an involution $\gt$ as in the
following diagram:
\[
\xymatrix{
& W \ar [dl] \ar [dr]^{/\gt} \\
C \ar [ddr]^f & & Z \ar [d]_{q_X} \ar[rd]_{q_F} \ar[rrd]^{q_Y}\\
& & X \ar [dl]^g \ar[dr]^{p_X} & F \ar[d]^{p_F} & Y, \ar[dl]^{p_Y}\\
&|K_C+\ga|^* & & E
}
\]
Eventually we will prove that the double cover $Y\ra E$ is the double cover from Theorem \ref{Tpc}.
\end{dsc}
\begin{dfn}\label{Dtrigo}
Since $h^0(K_C+\ga)=2$ the linear system $|K_C+\ga|$ is a line.
Denote the degree $4$ cover $C\ra|K_C+\ga|^*$ by $f$.
Define:
\[
  W:=\ol{C\times_{|K_C+\ga|^*}C\sm\mathrm{diagonal}},
\]
where both morphisms are $f:C\ra|K_C+\ga|^*$. The curve $W$ can be realized as the correspondence curve
\[
  W=\{(p_1,p_2)|p_1+p_2<K_C+\ga\}\subset C\times C.
\]
The curve $W$ admits an involution
\[
  \gt:=(p_1,p_2)\mapsto(p_2,p_1).
\]
Define
\[
  Z:=W/\gt.
\]
The curve $Z$ has a natural involution (which we later denote by $\gs$)
defined as follows:
Let $r$ be a point in $|K_C+\ga|^*$, and let $\{p_1,p_2\}\subset f^{-1}(r)$ be a degree $2$ divisor on the curve $C$.
The involution takes $\{p_1,p_2\}$ to $f^{-1}(r)\sm\{p_1,p_2\}$.
Denote by $X$ the quotient by this involution. By the trigonal
construction dictionary (see \cite{Do} p. 74), the double cover $Z\ra X$ is
admissible.
\end{dfn}
\begin{cor}
The arithmetic genera of the curves $Z,X$ are $7,4$ respectively.
\end{cor}
\begin{proof}
This follows from the trigonal construction dictionary (see \cite{Do} p. 74).
\end{proof}
\begin{ntt}\label{NHC}
If $C$ is hyperelliptic, we denote the hyperelliptic divisor by $H_C$.
\end{ntt}
\begin{lma}\label{LHKC}
If the curve $C$ is hyperelliptic, then the hyperelliptic divisor $H_C$ does
not satisfy:
\[
  H_C<K_C+\ga.
\]
\end{lma}
\begin{proof}
Assume that the above property holds. Since
\[
  h^0(H_C)=2=h^0(K_C+\ga),
\]
the linear system $|K_C+\ga|$ has base points (which sum to $H_C+\ga$). This is a
contradiction to the requirement that the linear system $|K_C+\ga|$ has no base
points.
\end{proof}
\begin{lma}\label{LisoSP}
The map
\[
  \begin{aligned}
  \Sym^2C&\ra\Pic^2C \\
  \{p_1,p_2\}&\mapsto p_1+p_2
  \end{aligned}
\]
is smooth on $Z$.
\end{lma}
\begin{proof}
The determinant of the Jacobian of this map vanishes exactly when $p_1+p_2$
is special. i.e. when
\[
  h^0(p_1+p_2)>1.
\]
In this case, the curve $C$ is hyperelliptic, and the hyperelliptic divisor $H_C$
satisfies:
\[
  H_C<K_C+\ga.
\]
By Lemma \ref{LHKC}, this does not happen in our setting.
\end{proof}
\begin{ntt}
Denote by $\gT$ the theta divisor of $\jac(C)$.
\end{ntt}
\begin{lma}\label{LZTheta}
The map
\[ 
  \begin{aligned}
  Z &     \ra \gT\cap(\gT+\ga) \\
  \{p_1,p_2\}&\mapsto p_1+p_2 \\
  \end{aligned}
\]
is an isomorphism.
\end{lma}
\begin{proof}
By the definition of $Z$, it
is a bijection on points. By Lemma \ref{LisoSP} it is an injection. To show
that it is an isomorphism, it is enough to calculate the genus of
$\gT\cap(\gT+\ga)$, which we do by adjunction:
\[
  \begin{aligned}
  K_{\gT}&=(\gT+K_{\Pic^2(C)})\gT=\gT^2 \\
  K_{\gT\cap(\gT+\ga)}&=(\gT\cap(\gT+\ga)+K_{\gT})(\gT\cap(\gT+\ga)) \\
    &\sim(\gT^2+\gT^2)\gT=2\gT^3=12.
  \end{aligned}
\]
Whence, the genus of $\gT\cap(\gT+\ga)$s is $\frac{12+2}{2}=7$.
\end{proof}
\begin{dfn}\label{Dij}
The maps
\[
  \begin{aligned}
  i:\gT\cap(\gT+\ga) &\ra    (\gT-\ga)\cap\gT=(\gT+\ga)\cap\gT \\
    d              &\mapsto d-\ga          \\
  j:\gT\cap(\gT+\ga) &\ra    (\gT-\ga)\cap\gT=(\gT+\ga)\cap\gT \\
    d              &\mapsto K_C-d          \\
  \gs:=i\circ j:\gT\cap(\gT+\ga) &\ra\gT\cap(\gT+\ga)\\
       d                       &\mapsto K_C+\ga-d,
  \end{aligned}
\]
are commuting involutions on $\gT\cap(\gT+\ga)$ (note that $\gs$ is the
quotient associated to the admissible double cover $Z\ra X$, see Definition
\ref{Dtrigo}).
\end{dfn}
\begin{rmr}
The involution $i$ has no fixed points.
\end{rmr}
\begin{ntt}\label{NE}
Denote
\[
  Y:=Z/i, \qquad F:=Z/j, \qquad E:=Z/\langle i,j\rangle.
\]
Denote the quotient maps from $Z$ to $E,X,Y,F$ by $q_E,q_X,q_Y,q_F$
respectively.
Since the involutions $i,j,\gs$ on the curve $Z$ commute (see Definition
\ref{Dij}), they induce involutions on the curves $X,Y,F$. Denote the
quotients by these involutions by $p_X,p_Y,p_F$
respectively.
Denote the involution induced by $j$ on $X$ by $j_X$.
Denote by $\gS$ the fixed points of the involution $j$ on the curve $Z$.
\end{ntt}
\begin{lma}\label{Lthetachar}
The theta characteristics of $C$ which are contained in
$Z=\gT\cap(\gT+\ga)$ are exactly the points of $\gS$. Moreover:
\begin{enumerate}
\item They are all odd.
\item They are paired, by the involution $i$, into pairs of the form $\{z,z+\ga\}$.
\item There are 12 such points.
\end{enumerate}
\end{lma}
\begin{proof}
By the definition of the involution $j$ (see Definition \ref{Dij}), the fixed
points of $j$ are exactly the theta characteristics of the curve $C$ in 
the curve $Z$.
By Definition \ref{Dij} these theta characteristics are paired by the
involution $i$.
Since the genus of the curve $C$ is 3, if there exists a theta characteristic
$\ga'$ such
that $h^0(\ga')=2$, then the curve $C$ is hyperelliptic, and $\ga'=H_C$, the
hyperelliptic divisor.
However, in this case $H_C<K_C+\ga$, and by Lemma \ref{LHKC} this is impossible.
Whence, all the theta characteristics in $\gS$ are odd (since they are
effective). Recall that the cover
\[
  \begin{aligned}
  \substack{\text{pairs of distinct}\\\text{odd theta characteristics}}&\ra\jac(C)[2]\sm\{0\} \\
  \{\theta_1,\theta_2\}&\mapsto \theta_1-\theta_2
  \end{aligned}
\]
is $SP_6(2)$ equivariant, and that the group $SP_6(2)$ acts transitively on
$\jac(C)[2]\sm\{0\}$. Whence, all the fibers of this map are of the same size:
$\left(\frac{28\cdot 27}{2}\right)/63=6$. By their definition, the pairs in
each fiber do not intersect.
\end{proof}
\begin{lma}\label{LZsing}
The singular points of the curve $Z$ are not in the set $\gS$. Moreover,
these points are paired by the involution $\gs$.
\end{lma}
\begin{proof}
By the trigonal construction dictionary (see \cite{Do} p. 74), the singular points of $Z=\gT\cap(\gT+\ga)$ are halves of $K_C+\ga$ in $Sym^2(C)$. However, in Lemma \ref{Lthetachar}
we proved that in the set $\gS$ there are only halves of $K_C$. By the trigonal
construction dictionary
singular points in $Z$ arise only if there are two points $p_1,p_2\in C$
such that $2(p_1+p_2)=K_C+\ga$. In this case
\[
  2(K_C+\ga-(p_1+p_2))=2K_C+2\ga-K_C-\ga=K_C+\ga.
\]
\end{proof}
\begin{ntt}\label{Nqis}
Denote the 6 pairs of odd half theta-characteristics in $\gS$ by 
$\{l_{i1},l_{i2}\}_{1\leq i\leq 6}$. Denote by $l^*_{ij}$ the images
in $|K_C|={\BP^2}^*$ of the $l_{ij}$.
Denote by $q_1\ldots q_6$ the points of $\gS/i$. The $q_i$s are
identified as the 6 ramification points of $Y\ra E$ over smooth points.
\end{ntt}
\begin{cor}\label{CFEadmis}
The double cover $F\ra E$ is admissible. The cover $Y\ra E$
is admissible over the nodes of the curve $E$.
\end{cor}
\begin{proof}
By Definition \ref{Dij}, the curve $Z$ can be realized as $F\times_E X$.
As the singular points of the curve $Z$ are paired by the involution $\gs$,
the admissibility of the double cover $F\ra E$ follows from the admissibility of the double cover $Z\ra X$.
By Lemma \ref{LZsing} the double cover $X\ra E$ is etale over nodes. Since
the double cover $F\ra E$
is admissible over nodes and since the curves $Y,F,X$ are the
intermediate quotients of the cover $Z\ra E$, the cover $Y\ra E$ is
admissible over nodes.
\end{proof}
\begin{cor}\label{Cgenera}
The arithmetic genera of the curves $E,F,Y$ are $1,1,4$ respectively.
\end{cor}
\begin{proof}
This follows by the Riemann-Hurwitz and counting singular points.
\end{proof}
\begin{dfn}\label{DCP2}
Denote by $k:C\hra|K_C|^*$ be the canonical embedding. Define the maps $\psi$
and $i_F$ by:
\[
  \begin{aligned}
  \psi:\Sym^2 C&\ra |K_C| \\
       \{p_1,p_2\}&\mapsto \ol{k(p_1)k(p_2)},\\
  i_F:F&\ra|K_C|\\
  \{\{p_1,p_2\},\{p_3,p_4\}\}&\mapsto\{p_1,p_2,p_3,p_4\}.
  \end{aligned}
\]
Note that by the definition of the curve $F$, the map
$\psi$ restricted to the curve $Z$ factors through the map $i_F$.
\end{dfn}
\begin{lma}\label{LFembed}
If the curve $F$ is not the irreducible nodal curve of arithmetic genus 1,
then the map $i_F:F\ra|K_C|$ is an embedding, and the image is a plane curve
of degree 3, with at most nodes.
\end{lma}
\begin{proof}
Through any point in $p\in C\subset|K_C|^*$ there are $3$ lines 
connecting it to each of the three points which complement $p$ to $K_C+\ga$.
Generically, these line are distinct.
Whence, the image of the map $i_F$ is a curve of degree $3$, and the map
$i_F$ is an immersion. A plane cubic which is not the cuspidal cubic has at
most nodes, and its arithmetic genus is 1. Since the arithmetic genus
of $F$ is 1, if the image of $i_F$ is
not the cuspidal curve then $i_F$ is an embedding. If the image is the
cuspidal curve then, since the arithmetic genus of the $F$ is $1$, the
curve $F$ is an  irreducible nodal curve of arithmetic genus $1$.
\end{proof}

%
\section{The Prym of the double cover}\label{Sprym}
%
\begin{dsc}
The section has two objectives. The first is to prove
(in proposition \ref{Psing}) that if the curve $E$ is irreducible, then it
is smooth.
The second objective is to calculate the variety 
$\prym(Y/E)$ in terms of the variety $\jac(C)$ and $\ga$, this objective is not essential
to the proof of Theorem \ref{Tbitan}.
\end{dsc}
\begin{ntt}
Denote by $\ti{Z},\ti{X},\ti{F},\ti{Y},\ti{E}$ the normalization of $Z,X,F,Y,E$
(see Definition \ref{Dtrigo} and Notation \ref{NE}) respectively.
\end{ntt}
\begin{lma}\label{Lnormisog}
The Norm map
\[
  \norm:\prym(\ti{Z}/\ti{X})\ra \prym(\ti{Y}/\ti{E})
\]
is an isogeny.
\end{lma}
\begin{proof}
We prove that the induced map of tangent spaces is an isomorphism.
Consider the homology groups $H^0(-,\gO^1_-)$ of the curves in the diagram:
\[
\xymatrix{
& \ti{Z} \ar[dr] \ar[d] \ar[dl] \\
\ti{X} \ar[dr] & \ti{F} \ar[d] & \ti{Y} \ar[dl] \\
& \ti{E}
}
\]
as $G=\gal(Z/E)$ modules.
Denote 
\[
  M:=H^0(\ti{Z},\gO^1_{\ti{Z}}).
\]
The group $G$ has only 3 non-trivial characters: $\chi_i,\chi_j,\chi_\gs$
Where
\[
  \ker(\chi_x)=\{1,x\}\text{ for all }x\in G\sm\{1\}.
\]
Denote by $\chi_1$ the trivial character and by $M_x$ the $\chi_x$
eigenspace of $M$. For any $G$-module $A$, denote by $A^x$ the elements of
$A$ fixed by the action of $x$.
The following identities hold:
\[
  \begin{aligned}
    H^0(\ti{Z},\gO^1_{\ti{Z}})&=M_1 \oplus M_i \oplus M_j\oplus M_\gs \\
    H^0(\ti{F},\gO^1_{\ti{F}})&=M_1^j\oplus M_i^j\oplus M_j^j\oplus M_\gs^j=
      M_1\oplus M_j,\\
    H^0(\ti{E},\gO^1_{\ti{Z}})&=M_1.
 \end{aligned}
\]
Since the genera of $\ti{E},\ti{F}$ are equal, then $\norm$ map
\[
  \norm:H^0(\ti{F},\gO^1_{\ti{F}})\ra H^0(\ti{E},\gO^1_{\ti{E}})
\]
is an isomorphism.
Whence, $M_j=0$. Moreover, we have:
\[
  \begin{aligned}
    H^0(\ti{Y},\gO^1_{\ti{Y}})&=M_1^i \oplus M_i^i \oplus M_\gs^i=M_1 \oplus M_i, \\
    H^0(\ti{X},\gO^1_{\ti{X}})&=M_1^{\gs}\oplus M_i^{\gs}\oplus M_\gs^{\gs}
      =M_1\oplus M_{\gs}.
  \end{aligned}
\]
Whence, $\norm:H^0(\ti{Z},\gO^1_{\ti{Z}})\ra H^0(\ti{Y},\gO^1_{\ti{Y}})$ induce an isomorphism
\[
  \begin{aligned}
    M_i\cong\ker&(\norm(H^0(\ti{Z},\gO^1_{\ti{Z}})\ra H^0(\ti{X},\gO^1_{\ti{X}}))) \\
 \cong\ker&(\norm(H^0(\ti{Y},\gO^1_{\ti{Y}})\ra H^0(\ti{E},\gO^1_{\ti{E}}))) .
  \end{aligned}
\]
Therefore it induces an isogeny on the $\prym$s.
\end{proof}
\begin{lma}\label{Lkernel}
The kernel of the norm map:
\[
  \norm:\prym(Z/X)\ra\prym(Y/E),
\]
is a subset of $\prym(Z/X)[2]$.
\end{lma}
\begin{proof}
Denote by $[2]_-$ the multiplication by $2$ operator on an Abelian group. Note
that:
\[
  {q_Y}_*q_Y^*:\jac(Y)\ra\jac(Y)=[2]_{\jac(Y)}.
\]
Denote by $\psi_*,\psi^*$ the restrictions of maps ${q_Y}_*,q_Y^*$ to
the Abelian varieties $\prym(Y/E),\prym(Z/X)$ respectively.
We then have
\[
  \begin{aligned}
  \psi_*\psi^*&=[2]_{\prym(Y/E)} &\Rightarrow \\
  (\psi^*\psi_*)^2&=\psi^*[2]_{\prym(Y/E)}\psi_*=\psi^*\psi_*[2]_{\prym(Z/X)}
  \qquad &\Rightarrow \\
  \psi^*\psi_*&=[2]_{\prym(Z/X)}.
  \end{aligned}
\]
Therefore the kernel of the norm map is a subset of $\prym(Z/X)[2]$.
\end{proof}
\begin{lma}\label{LWirr}
If the curve $W$ is reducible then the cover $C\ra|K_C+\ga|^*$ has an intermediate factor.
\end{lma}
\begin{proof}
If the curve $W$ is reducible then the Galois group of the Galois closure of
the cover $C\ra\BP^1$ is either 
of $D_4,\Cy{4},\Cy{2}\times\Cy{2}$. In each of these cases the cover
$C\ra|K_C+\ga|^*$ has an intermediate factor.
\end{proof}
\begin{prp}\label{Psing}
If the curve $E$ is irreducible, then it is smooth.
\end{prp}
\begin{proof}
Assume that the curve $E$ is singular and irreducible, Since its arithmetic
genus is 1, it has one node. By Corollary \ref{CFEadmis} The curves $F,Y$
are admissible covers of $E$. Thus each of them has one node.
By Lemma \ref{Lkernel} the kernel $I$ of the norm map
\[
  \jac(C)\cong\prym(Z/X)\ra\prym(Y/E)\cong\jac(\ti{Y})
\]
is a subgroup of $\jac(C)[2]$. 
The action of
$SP_6(2)$ on the odd theta characteristics is 2-transitive (see \cite{CCNPW},
p. 46, \cite{DM} 246-248). By the monodromy argument (Lemma \ref{Lmonodromy}),
the kernel $I$ of the isogeny above can only be one of the following subspaces:
\[
    0, \jac(\ti{Y})[2].
\]
In either case, we have
\[
  \jac(\ti{Y})\cong\jac(C).
\]
By Torelli, $C\cong\ti{Y}$. i.e. the curve $C$ is hyperelliptic.
Denote by $p_1,p_2$ the two points in $\ti{Y}$ which lie over the singular point of $Y$. We will use two consecutive monodromy arguments to show that
$\ga+H_C\sim p_1+p_2$:
Since the stabilizer of the odd theta characteristic $p_1+p_2$ in $SP_6(2)$ is
transitive on the even theta characteristics (see \cite{CCNPW} p. 46), the
theta characteristic $\ga+H_C$ is odd.
Since the group $SP_6(2)$ is 2-transitive on the odd theta characteristics
(see \cite{DM} p. 247) the equality $\ga+H_C\sim p_1+p_2$ holds.
Therefor the linear system $K_C+\ga$ has two base points.
\end{proof}
\begin{thm}\label{TZiEC}
If the curve $E$ is smooth, then under the identification $\prym(Z/X)\cong\jac(C)$ (see \cite{Do} Theorem 2.11),
the kernel of the norm map 
\[
  \jac(C)\cong\prym(Z/X)\ra\prym((Z/i)/E),
\]
is $\ga^\perp$.
\end{thm}
\begin{proof}
The identification $\prym(Z/X)\cong\jac(C)$, identifies
$\prym(Z/X)[2]\cong\jac(C)[2]$.
By monodromy considerations, it can only be one of the following subspaces:
\[
    0, \langle\ga\rangle, \ga^\perp, \jac(C)[2].
\]
The double cover $Y\ra E$ has exactly  6 ramification points.
Choose some partition to pairs of these
ramification points. The singularization on these pairs (see Notation
\ref{Nsing}) is an admissible cover. By Lemma \ref{LDLisogeny},
the polarization type of  $\prym(Y/E)$ is $2^1 1^2$.
The only possibility for the kernel is therefore the group $\ga^\perp$.
\end{proof}

%
\section{The plane configuration}\label{Splane}
%
\begin{dsc}
We have already identified two of the moduli spaces appearing in the statement
of Theorem \ref{Tpc}. In this section we identify the third moduli space, and
prove Theorem \ref{Tpc}. In Theorem \ref{TChee} we identify the
hyperelliptic locus (where the theta characteristic $\ga+H_C$ is even).
Throughout this section we assume that the cover $C\ra|K_C+\ga|^*$ does
not have an intermediate factor. Recall that by Lemma \ref{LWirr} and
Proposition \ref{Psing} the curve $E$ is smooth.
\end{dsc}
\begin{lma}\label{LYE11}
Let $R$ be a set of 6 distinct points on $E$. There is an
isomorphism between the coarse moduli spaces describing
the following sets of data:
\begin{itemize}
  \item Embeddings of $E$ in $\BP^2$ such that $R$ is contained in a conic $Q$;
    up to projective transformations.
  \item Double covers $\pi:S\ra E$, with an associated involution $i_\pi$, such
    that the ramification divisor is $R$; up to isomorphisms.
\end{itemize}
\end{lma}
\begin{proof}
\item {\em From $Q$ to $S$:} Let $K(E)$ be the function field of $E$.
Denote By $\eta_Q$ the function defining $Q$ in $\BP^2$, and by $\eta_l$ a linear
form on $\BP^2$. Let $S$ be the smooth projective model of
\[
  K(E)(\sqrt{\eta_Q/\eta_l^2}).
\]
\item {\em From $S$ to $Q$:} Consider the group $\{\Id,i_\pi\}$ of
automorphisms of $S$. In a similar fashion to the proof of Lemma
\ref{Lnormisog}, the characters of this group induce a decomposition of
$\pi_*\cO_S$. The $1$ eigenspace is $\cO_E$, so we have:
\[
  \pi_*\cO_S\cong\cO_E\oplus\cO_E(-L),
\]
for some divisor $L$ on $E$.
However, we now have (see \cite{Ha} p. 306)
\[
  \cO_E(-R)\cong\det(\pi_*\cO_S)=\cO_E(-L)^{\otimes 2}=\cO_E(-2L),
\]
Or $R\sim 2L$. This means that under the embedding given by the linear series
$|L|$:
\[
  \phi_{L}:E\ra\BP^2,
\]
the points of $R$ sit on a conic $Q$. 
These constructions are the inverse of one another.
\end{proof}
\begin{ntt}\label{NQS}
Given a double cover $\pi:S\ra E$ and an embedding $E\hra\BP^2$ as in
Lemma \ref{LYE11}, denote the $Q$
from Lemma
\ref{LYE11} by $Q_S$ (it is defined up to the projective transformations that
fix $E\subset\BP^2$ as a set). Denote the divisor $L$ from the proof of
Lemma \ref{LYE11} by $L_S$.
\end{ntt}
\begin{lma}\label{LYcan}
Given plane curves $Q_S,E$ as above, denote by $\eta_Q,\eta_E\in\BF[x_0,x_1,x_2]$
homogeneous polynomials for $Q_S,E$ respectively. The curve $S$ is the complete
intersection of the following surfaces in $\BP^3$. 
\begin{itemize}
  \item $Q_3$, a cubic given by
$\ti{\eta_E}:=\eta_E\in\BF[x_0,x_1,x_2]\subset\BF[x_0,x_1,x_2,x_3]$
  \item $Q_2$, a quadric given by $\ti{\eta_Q}:=\eta_Q+x_3^2\in\BF[x_0,x_1,x_2,x_3]$
\end{itemize}
Moreover, this intersection is transversal.
\end{lma}
\begin{proof}
By the proof of Lemma \ref{LYE11} The curve $S$ is the complete intersection
$Q_2\cap Q_3$. Since One of the surfaces is a quadric and the other is a cubic,
The complete intersection $Q_2\cap Q_3$ is the canonical embedding of $S$
in $\BP^3$.
We will check that the surfaces $Q_2,Q_3$ intersect transversely by
calculating the rank of the Jacobian matrix of the intersection:
\[
  \left(\begin{array}{c}
   \nabla \ti{\eta_Q} \\
   \nabla  \ti{\eta_E}
  \end{array}\right)=
  \left(\begin{array}{c|c}
   \nabla\eta_Q  & 2x_3 \\
   \nabla\eta_E  & 0 
  \end{array}\right).
\]
\item{\em The case $x_3\neq 0$:} Since $E$ is smooth, $\nabla\eta_E\neq 0$ and:
\[
  \mathrm{rank}
  \left(\begin{array}{c|c}
   \nabla\eta_Q& 2x_3 \\
   \nabla\eta_E& 0 
  \end{array}\right)
=1+\mathrm{rank}(\nabla \eta_E)=2
\]
\item{\em The case $x_3=0$:} Since $E,Q_S$ intersect transversally,
\[
  \mathrm{rank}
  \left(\begin{array}{c}
   \nabla \eta_Q \\ \nabla\eta_E 
  \end{array}\right)=2.
\]
\end{proof}
\begin{ntt}
Given plane curves $Q_S,E$ as above denote by $Q_2(S/E),Q_3(S/E)$ the
surfaces $Q_2,Q_3$ appearing in the proof of Lemma \ref{LYcan}.
\end{ntt}
\begin{ntt}\label{Ngb}
Recall (Lemma \ref{Cgenera}) that the genera of the curves $E,F$ are both 1.
Denote by $\gb$ generator of $\ker(F\ra E)$, and by $\gb'$ the unique generator
of the image of $\Pic(F)[2]$ under the double cover $F\ra E$.
\end{ntt}
\begin{lma}\label{LXY}
There is a natural isomorphism of sheaves over the curve $E$:
\[
  {p_Y}_*(\cO_Y)\otimes\cO_E(\gb')={p_X}_*(\cO_X).
\]
\end{lma}
\begin{proof}
By Definition \ref{Dij}, the curve $Z$ can be realized as the fibered product
$Y\times_E F$. Whence (recall Notation \ref{NQS}):
\[
  \begin{aligned}
   (q_E)_*\cO_Z=&{p_Y}_*\cO_Y\otimes{p_F}_*\cO_F \\
   =&(\cO_E\oplus\cO_E(-L_Y))\otimes(\cO_E\oplus\cO_E(\gb')) \\
   =&\cO_E\oplus\cO_E(\gb')\oplus\cO_E(-L_Y)\oplus\cO_E(-L_Y+\gb').
  \end{aligned}
\]
Decomposing this sheaf to irreducible representations of $Gal(Z/E)$ (as in
the proof of Lemma \ref{Lnormisog}), the result follows.
\end{proof}
\begin{cor}\label{CXY}
The following identity hold:
\[
  L_Y-L_X=\gb'.
\]
\end{cor}
\begin{proof}
This follows from Lemma \ref{LXY}.
\end{proof}
\begin{proof}[Proof of Theorem \ref{Tpc}]
By the results of this section, we see that there is an equivalence between
the following sets of data
\begin{itemize}
  \item Ramified double covers $S/E$, of irreducible smooth curves of
    genera $4,1$ respectively, such that $E$ is smooth.
  \item Configurations of plane curves $Q_S,E$ of degrees 2,3 respectively,
    where $E$ is smooth, and $E,Q_S$ intersect transversally.
  \item configurations of $Q_2(S/E),Q_3(S/E)$: a quadric and a cubic in
    $\BP^3$, 
    which arise from a plane configuration.
\end{itemize}
The curve $C$ defines the pair $((Y\ra E),\gb')$ as in Section \ref{SCZ}.
By Proposition \ref{Psing}, the curve $E$ is smooth. By the trigonal
construction dictionary (see \cite{Do} p. 74) the curve $X$ is smooth.
By the proof of Lemma \ref{LYcan} the intersection $Q_X\cap E$ is transversal.
By Corollary \ref{CXY}
the intersection $Q_Y\cap E$ is transversal, and by the proof of
Lemma \ref{LYcan} the
curve $Y$ is smooth.
To show that this map is a bijection (on the coarse moduli spaces) 
we construct the inverse of the map:
observe that
\begin{itemize}
\item The curve $Z$ satisfies $Z=Y\times_E(E/\gb')$.
\item By Corollary \ref{CXY}, the curve $X$ is uniquely defined by the
  double cover $Y\ra E$ and the element$\gb'\in\Pic(E)[2]\sm\{0\}$.
\item Generically, the curve $X$ has exactly two $g^1_3$s, which are conjugate
under the action of $j_X$ (recall the definition of the involution $j_X$ in
Notation \ref{NE}). Therefor these $g^1_3$s are given as the rulings of
the quadric $Q_2(X/E)$. Taking flat limits the map $g: X\ra\BP^1$ is always
one of the $g^1_3$s appearing as a ruling of the quadric
$Q_2(X/E)$. There are at most two such $g^1_3$s, and they are conjugate under
the involution $j_X$.
\item The map $f:C\ra\BP^1$ is the trigonal construction on the tower
\[
  Z\ra X\overset{g}{\ra}\BP^1.
\]
\end{itemize}
\end{proof}
\begin{prp}\label{Pbase}
A pair $(C,\ga)$ fails the requirements of Notation \ref{NCa} if and only if
the curve $C$ is hyperelliptic and the theta characteristic $\ga+H_C$ is odd.
\end{prp}
\begin{proof}
\item {\em The non-hyperelliptic case:} For different indices $i$, the bitangents
$l^*_{i1},l^*_{i2}$ (see Notation \ref{Nqis}) intersect the canonical model of the curve $C$ at different points. By the definition of $l_{ij}$, the four bitangency points of
$l_{i1}\cup l_{i2}$ to $C$ sum to $K_C+\ga$.
\item {\em The hyperelliptic, even $\ga+H_C$ case:} Since $\ga+H_C$ is even there
are two non-intersecting sets of four Weierstrass points (each) which sum to
$K+\ga$.
\item {\em The hyperelliptic, odd $\ga+H_C$ case:} By assumption, there exist
two points $p_1,p_2\in C$ such that $\ga+H_C=p_1+p_2$. Whence, $K_C+\ga=H_C+p_1+p_2$.
\end{proof}
\begin{thm}\label{TChee}
The following statements are equivalent:
\begin{itemize}
  \item The curve $C$ is hyperelliptic and $\ga+H_C$ is an even theta
    characteristic.
  \item the conic $Q_X$ (see Notation \ref{NQS}) is singular.
\end{itemize}
\end{thm}
\begin{dsc}
The proof of this theorem is not direct. We will use some preliminary lemmas:
\end{dsc}
\begin{lma}\label{LQsing}
If $Q_X$ is singular, then there exists a map $t_1:E\ra\BP^1$ of degree 3,
and a map $t_2:\BP^1\ra\BP^1$ of degree 2 such that:
\[
  t_1\circ p_X=t_2\circ g.
\]
(Recall the definitions of $g,p_X$ in \ref{trigo} and Notation \ref{NE}
respectively)
\end{lma}
\begin{proof}
Since the surface $Q_X$ is singular, the quadric $Q_2(X)$ has only one ruling.
Therefore, with the appropriate choice of coordinates (see Notation \ref{NE}),
\[
  g=j_X(g)\qquad \Rightarrow \qquad j_X(g)\cdot g=g^2.
\]
\end{proof}
\begin{ntt}\label{NGprime}
We consider the group $S_4\times(\BZ/2)$. We denote by
$p_1$ and $p_2$ the projections on the first and second factors.
For any subgroup $G\subset S_4$, denote
\[
  G':=\{(\gc,\text{sign}(\gc))|\gc\in G\}.
\]
\end{ntt}
\begin{lma}\label{LS4Z2}
Let $G$ be a subgroup of $S_4\times(\BZ/2)$ such that
\[
   G\not\supset\{1\}\times\Cy{2}, \qquad
   G\not\subset S_4\times\{1\}, \qquad
   p_1(G)\supset S_3,
\]
Then $G=S_3'$.
\end{lma}
\begin{proof}
If $G$ is a subgroup of $S_4$, there are exactly 3 subgroups
$\ti{G}\subset S_4\times(\BZ/2)$ such that $p_1\ti{G}=G$. These are:
\[
  G\times\{1\}, \qquad G\times\Cy{2}, \qquad G'.
\]
\end{proof}
\begin{lma}\label{L5p}
Let $\cM\subset\cM_3$ be a 5-dimensional subscheme that  parameterize a family
of smooth curves of genus 3 with involutions,
then $\cM$ is a subset of the hyperelliptic locus of $\cM_3$.
\end{lma}
\begin{proof}
The quotients by an involution on a curve of genus 3 can be of genus 0,1 or
2. By Riemann-Hurwitz, in these cases, there are 8,4 and 0 branch points
respectively. Therefore, there are only 4 parameters if the genus of the
quotient is 1,
and only 3 parameters if the genus of the quotient is 2. Since a smooth limit
of a hyperelliptic curve is hyperelliptic, the assertion is proved.
\end{proof}
\begin{proof}[Proof of Theorem \ref{TChee}]
Since $Q_X$ is singular, the quadric $Q_2(X/E)$
has only one ruling (recall form Lemma \ref{LWirr} that $X$ is irreducible).
By Lemma \ref{LQsing}, the map $g\cdot j_X(g)$ factors through
$E$. The situation can be summarized in the following commutative diagram:
\[
  \xymatrix{
          & W \ar[ddl] \ar[dr]^{/\gt} \\
          &                   & Z \ar[d]^{q_X} \ar[dr]^{q_F} \\
C \ar[dr]^f&                   & X \ar[dl]^g \ar[d]^{p_X} & F. \ar[dl]^{p_F} \\
          & \BP^1_{hi} \ar[d]^{t_1}              & E \ar[dl]^{t_2} \\
          & \BP^1_{lo}
  }
\]
Denote by $\ol{F},\ol{W}$ the Galois closure of the covers
$F/\BP^1_{lo},W/\BP^1_{lo}$ respectively.
\item{\em Claim 1. The Galois group of $\ol{F}/\BP^1_{lo}$ is a subgroup of $S_4$:}
Since $F/E$ is unramified, this follows from the trigonal construction
properties (see \cite{Do} p. 71-72).
\item{\em Claim 2. The equality $\ol{W}=\ol{F}\times_{\BP^1_{lo}}\BP^1_{hi}$ holds:}
Since $X=E\times_{\BP^1_{lo}}\BP^1_{hi}$, there exists a cover
\[
  \ol{F}\times_{\BP^1_{lo}}\BP^1_{hi}\ra Z.
\]
Since the Galois closure of $Z/\BP^1_{hi}$ is the Galois closure of
$C/\BP^1_{hi}$, there exists a cover
\[
  \ol{F}\times_{\BP^1_{lo}}\BP^1_{hi}\ra C.
\]
Hence, there exists a cover
\[
  \ol{F}\times_{\BP^1_{lo}}\BP^1_{hi}\ra W.
\]
This gives a bound from above on the Galois closure. However, since the Galois
group of the Galois closure of $C/\BP^1_{hi}$ is $S_4$, we know that
\[
  Gal(\ol{W}/\BP^1_{lo})\supsetneq S_4.
\]
\item{\em Claim 3. The group corresponding to $C$ in the Galois extension
$\ol{W}/\BP^1_{lo}$ is $S_3'$ (see Notation \ref{NGprime}):}
This group satisfies the requirements of Lemma \ref{LS4Z2}.
\item{\em Claim 4. The curve $C$ is hyperelliptic, and $\ga$ is even:}
Since
\[
  [S_3\times\Cy{2}:S_3']=2,
\]
the curve $C$ is a double cover.
Since the curve $X$ has 5 parameters, the curve $C$ has 5 parameters too.
By Lemma \ref{L5p}, the curve $C$ is hyperelliptic.
By Proposition \ref{Psing}, the theta characteristic $\ga$ is even.
\item{\em Claim 5: For any hyperelliptic $C$ and $\ga$ even, The quartic $Q_X$ is singular:} 
This follows from irreducibility and dimension arguments on the moduli space
\[
  \{(C,\ga)|\text{genus}(C)=3,\ga\in\jac(C)[2]\sm\{0\},\ga+H_C\text{ is even}\}.
\]
\end{proof}
%
\section{Bitangents of a canonically embedded quartic}\label{Strigo}
%
\begin{dsc}
In this section we assume that the curve $C$ is not hyperelliptic, and
that the cover $C\ra|K_C+\ga|^*$ does not have an intermediate factor. By
Proposition
\ref{Pbase}, for any $\ga$, the pair $(C,\ga)$ fulfills the requirements of
Notation \ref{NCa}.
Recall (see Proposition \ref{Psing} and Theorem \ref{TChee}) that both curves
$E,Q_X$ are smooth and intersect transversally. Our main goal in
is section is Theorem \ref{TP2prop} in which we describe
the plane data $(E,Q_Y\subset|\BP^2),\gb'\in\Pic(E)[2]$ in terms of the
$l_{ij}$s.
\end{dsc}
\begin{dfn}[The chord construction (c.c.)]\label{Dcc}
Let $V$  be a smooth curve of genus 1, and
$\gc\in\Pic(V)[2]\sm\{0\}$. Denote
$V':=V/\gc$. Let $\xi$ be the quotient map $V\ra V'$. Denote
by $\gc'$ the unique non-zero element in $\xi_*(\Pic(V)[2])$
(we then have $V'/\gc'\cong V$).
Let $\mu_V:V\ra\BP^2$ be a smooth embedding.
Define
\[
  \begin{aligned}
  \mu_{V'}:V'  &\ra{\BP^2}^* \\
   \{p,p+\gc\}&\mapsto\ol{\mu_V(p)\mu_V(p+\gc)}
  \end{aligned}
\]
\end{dfn}
\begin{lma}\label{Lcc}
The $c.c.$ is a smooth embedding, moreover, if $\mu_V(p)$ is a flex, then so is $\mu_{V'}(\pi(p)+\gc')$.
\end{lma}
\begin{proof}
This is a classical theorem. See \cite{LL} Theorems 3 and 4, or \cite{Sa} chapter V section V or \cite{Dol2}.
\end{proof}
\begin{prp}\label{PP2prop}
In the notations of Definition \ref{DCP2} the following properties hold:
\begin{enumerate}
 \item The involution induced by $i$ on the points of the curve $F$ acts in the
     following way: Let $\{p_1,p_2\}$ be a point in $Z$, and let $\ol{p_1 p_2}$
     be the corresponding line in the curve $F\subset|K_C|$.
    There is a unique pair of points
     $\{p_3,p_4\}$ s.t. $K_C+\ga\sim p_1+p_2+p_3+p_4$. We then have 
\[
  i(\ol{p_1 p_2})=\ol{p_3 p_4}.
\]  
  \item The curve $E$ embeds naturally in $|K_C|^*$ as the c.c. of $i_F$, (as
    the intersection point of each of the pairs of lines $\{l,i(l)\}$).
  \item The ramification points of $Z\ra F$ in $|K_C|\cong{\BP^2}^*$ are
    the 12 bitangents to $C$, which match the 12 effective theta
    characteristics in $Z$. (see lemma \ref{Lthetachar})
  \item The ramification points of the double cover $X\ra E$ are in 1-1
   correspondence with the pairs $\{l,i(l)\}$ of bitangents to $C$ in $F$.
\end{enumerate}
\end{prp}
\begin{ntt}\label{NsP1}
Denote by $s_{\BP^1}$ the map:
\[
  \begin{aligned}
  s_{\BP^1}:\BP^1\times\BP^1&\ra\Sym^2\BP^1\cong\BP^2 \\
    ((x_0;x_1),(y_0;y_1))&\mapsto(x_0 y_0;x_0 y_1+x_1 y_0;x_1 y_1).
  \end{aligned}
\]
\end{ntt}
\begin{lma}\label{LEembed}
The quotient $p_X:X\ra E$ can be extended to a quotient
map $Q_2(X/E)\ra\BP^2$. Moreover, in suitable coordinates on
$\BP^3\supset Q_2(X/E)$, this quotient map can be identified with the
quotient map $s_{\BP^1}$ (from notation \ref{NsP1}).
\end{lma}
\begin{proof}
\item{\em Claim 1. Existence of the extension of the quotient map on $Q_2(X/E)$:}
Keeping the coordinates of the proof of Lemma \ref{LYcan}, $\BP^3$ is realized
as a cone over $\BP^2$. By the proof of Lemma \ref{LYcan}, the vertex of this
cone is not in the surface $Q_2(X/E)$. Whence, each chord of the cone
intersects $Q_2(X/E)$ at two non-vertex points (counting multiplicities).
The projection of the base of the cone ($\BP^2$) is a quotient which
restricts to the involution $j_X$ on the curve $X$.
\item{\em Claim 2. Identification of the quotient maps $Q_2(X/E)\ra\BP^2$ and $s_{\BP^1}$:}
Choose coordinates on $\BP^3$ such that $Q_2(X/E)$ is given by the zeros of
$x_0 x_3=x_1 x_2$.
The quadric $Q_2(X/E)$ is the image of
\[
  \begin{aligned}
  \BP^1\times\BP^1 &\ra \BP^3 \\
  ((u_0;u_1),(v_0;v_1))&\mapsto (u_0 v_0;u_0 v_1;u_1 v_0;u_1 v_1).
  \end{aligned}
\]
The involution $s_{\BP^1}$ identifies
$(x_0;x_1;x_2;x_3)$ with $(x_0;x_2;x_1;x_3)$.
Fixing $\BP^2:=Z(x_1-x_2)$ and infinity as $(0;1;-1;0)$,
we see that the property holds.
\end{proof}
\begin{thm}\label{TXZi}
The ramification points of the double cover $X\ra E$ sit under the embedding
$E\subset\BP^2$ from Proposition
\ref{PP2prop} on the conic $Q_Y$ (see Notation \ref{NQS}).
\end{thm}
\begin{proof}
Given a canonical map $k:C\hra\BP^2$, The following commutative diagram 
describes the construction as presented in Proposition \ref{PP2prop}:
\[
  \xymatrix{
\BP^1\times\BP^1\cong Q_2(X/E) \ar[d]^{s_{\BP^1}} & X \ar[l]_<<<{g,g\circ j_X} \ar[d] &
Z \ar[d]^{q_F} \ar[l]_{q_X} \ar[r]^{i_Z} & \Sym^2 C \ar[d]^{\psi} \\
\Sym^2\BP^1\cong\BP^2 & E \ar[l] & F \ar[l]_{p_F} \ar[r]^{i_F} &
{\BP^2}^*=|K_C|
}
\]
(Note that $Q_2(X/E)\cong\BP^1\times\BP^1$ by Theorem \ref{TChee}).
We identify two embeddings $E\ra\BP^2$. The first is $i_E$ in the diagram.
The second is the c.c. (see Definition \ref{Dcc})of $i_F$ and $\gb$.
Let us fix some $q\in C$ and consider the degree 3 divisor on $Z$:
\[
  D:=\{\{p_1,p_2\}|\{p_1,p_2\}\in Z\}\subset Z.
\]
The images of $D$ in ${\BP^2}^*$ and $\BP^2$ are subsets of the following
hyperplanes (respectively).
\[
  \begin{aligned}
  \psi(i_Z(D))&\subset\{\{\ol{k(p_1)k(p_2)}\}|p_1\in C\}=k(p_2)^* \\
  s_{\BP^1}(\phi(D))&= s_{\BP^1}(\{(g(p_1+p_2),g(K_C-p_1-p_2))|\{p_1,p_2\}\in Z\})\\
                    & \subset s_{\BP^1}(\{f(p_2)\}\times\BP^1)
  \end{aligned}
\]
By the definition of the quotient $s_{\BP^1}$ (see Notation \ref{NsP1}), the
last divisor above is a line in $\BP^2$.
Let $l$ be a line in $|K_C|$. As proved above, the image of the set $F\cap l$
under the map $i_E\circ\pi$ are three collinear points in $\BP^2$. The divisor
class of the sum of these 3 points is $L_X$.
Denote by $l'$ a line in $|K_C|^*$, and by $\ti{E}$ the image of the c.c of
$i_F$ and $\gb$. By Lemma \ref{Lcc}:
\[
 L_X+\gb'=l'\cap\ti{E}.
\]
However, by Corollary \ref{CXY}, $L_Y=L_X+\gb$. Since the ramification locus
of $s_{\BP^1}$ is $Q_X$, we are done.
\end{proof}
\begin{thm}\label{TP2prop}
For pairs $(C,\ga)$ such that the curve $C$ is non-hyperelliptic and such that
the cover $C\ra|K_C+\ga|^*$ does not have an intermediate factor, the following
properties hold:
\begin{enumerate}
\item The curve $F$ is naturally isomorphic to the unique degree 3 curve
      passing through the points $l_{ij}^*$.
\item There exists $\gb\in\Pic(F)[2]$, s.t. $l^*_{i1}-l^*_{i2}=\gb$ for
     all $1\leq i\leq 6$.
\item The curve $E$ is naturally isomorphic to the chord construction on the
   curve $F$ and $\gb$.
   Under this construction, the points $l_{ij}^*$ collapse to 6 points on a
   conic.
\item The information encoded in each of the following objects is the same:
\begin{itemize}
  \item The 12-tuple $\{l^*_{ij}\}_{1\leq i\leq 6,j=1,2}$ 
    (viewed as points in $|K_C|$), with the choice of a pairing, 
    both determined by $\ga$.
  \item The curves $E,Q$ embedded in $\BP^2$, and a choice of an element
   $\gb'\in\Pic(E)[2]$
\end{itemize}
\end{enumerate}
\end{thm}
\begin{proof}
Recall that by Proposition \ref{PP2prop} the following properties hold:
\begin{itemize}
  \item The 12-tuple $\{l_{ij}^*\}$ is contained in the degree $3$ smooth curve $F$.
  \item The pairing on the 12-tuple $\{l_{ij}^*\}$ defines an element
\[
  \gb\in\Pic(F)[2]\sm\{0\}.
\]
  \item The image of the c.c (see Definition \ref{Dcc}) applied to the curve
$F$, and the divisor $\gb$ is the curve $E$.
\end{itemize}
In Theorem \ref{TXZi} we proved that the c.c on the curve
$F$, and the divisor $\gb$ sends the 12-tuple $\{l_{ij}^*\}$ to 6 points
on the conic $Q_Y$.
The double cover $F\ra E$ determines a unique divisor
\[
  \gb'\in\Pic(E)[2]\sm\{0\}
\]
(see Notation \ref{Ngb}). By Theorem \ref{Tpc}, the
information in $E,Q_Y,\gb'$ is in 1-1 correspondence with the
information in $(C,\ga)$.
\end{proof}
\begin{prp}\label{P9pts}
Any 9 of the bitangents $l_{ij}$, with at least one marked pair,
determine the other 3 bitangents, and the pairing on all 12 bitangents.
\end{prp}
\begin{proof}
The 9 points $l_{ij}^*$ define a unique cubic.
The marked pair determines a chord construction. By Theorem \ref{TXZi},
the images
of the 12 points $l_{ij}^*$ under the c.c. sit on the conic $Q_Y$.
However, as 5 points determine a
conic, and the c.c. is a double cover, the conic $Q_Y$ is defined by 
the images of the 9 given points. There are 6 points in the set $E\cap Q_Y$,
which pull back to 12 points under the double cover $F\ra E$.
\end{proof}
\begin{prp}\label{Pauto}
There is no non-trivial projective transformation that fixes the set of
bitangents $l_{ij}$ (as a set).
\end{prp}
\begin{proof}
Such a transformation stabilizes $F$, and the dual transformation stabilizes
both $E$ and $Q$.
\end{proof}
%
\section{Curves with automorphisms}\label{Ssing}
%
\begin{dsc}
In this section we analyze the degeneration of Theorem \ref{Tpc}
to the case where the cover $C\ra|K_C+\ga|^*$ has an intermediate factor.
Recall that by Lemma \ref{LWirr}, either the cover
$f:C\ra\|K_C+\ga|^*$ is Galois, or the Galois group of the Galois closure is
the group $D_4$.
Recall also that by Lemma \ref{LFembed}, the map $F\ra|K_C|$ is an embedding,
and the image is a curve of degree 3.
\end{dsc}
\begin{prp}
If the cover $f:C\ra|K_C+\ga|^*$ is Galois then
\begin{enumerate}
\item The curve $F$ is a union of $3$ lines, mapped under $i_F:F\ra|K_C|$ to
a union of 3 lines.
\item There are exactly $4$ bitangents on each of the components
\item The pair $C,\ga$ can be reconstructed from the curve $F$ and the 12
marked points.
\end{enumerate}
\end{prp}
\begin{proof}\label{Pgal}
By Galois correspondence,
\[
  W=C\coprod C\coprod C.
\]
Whence $F$ is a union of 3 lines. By symmetry considerations, there are
exactly 4 bitangents on each of them.
\end{proof}
\begin{prp}\label{PD4}
If the cover $f:C\ra|K_C+\ga|^*$ is not Galois, and the Galois group of the
Galois Closure is not $S_4$ then
\begin{enumerate}
\item The Galois group is $D_4$.
\item The curve $F$ is a union of two lines, mapped under $i_F:F\ra|K_C|$
to a union of a conic and a line.
\item There are $8$ bitangents on the component which maps to a conic,
and $4$ bitangents on the other component.
\end{enumerate}
\end{prp}
\begin{proof}
By Galois correspondence the curve $C$ has a double quotient. Denote this
double quotient by $C'$.
Define the tower of curves:
\[
  Z_0\ra X_0\ra\BP^1
\]
to be the bigonal construction on the tower
\[
  C\ra C'\ra\BP^1.
\]
By Galois correspondence, the curves $Z,X$ decompose to
\[
  Z=Z_0\coprod C',\qquad X=X_0\coprod \BP^1.
\]
Therefor the curve $E$ is a union of two lines. Since the cover $F\ra E$
is admissible (see \ref{CFEadmis}), the curve $F$ is a union of two lines.
Since the case describe in proposition \ref{Pgal} is a degeneration of the
current case, there are exactly $4$ bitangents on one of the
components. Whence, there are $8$ bitangents on the other components.
As the arithmetic genus of $F$ is $1$ the map $i_F$ maps one of these
components
to a line, and the other component to a conic. Since component with $8$
$l_{ij}$s
degenerate in proposition \ref{Pgal} to a conic, it is a conic.
\end{proof}
\begin{lma}\label{Lsingrecon}
The curve $C\subset|K_C|^*$ is determined by the set $l_{ij}$, with the pairs
marked on it.
\end{lma}
\begin{proof}
The curve $F$ is the unique cubic through the $l_{ij}$s. Since the
irreducible components of the curve $F$ are rational, the curve $Z$ is
the unique double cover of the curve $F$ with the given ramification points
$l_{ij}$. Applying the same argument on the curve $E$ the 
curve $X$ is also well defined. Thus the curve $X$ is the unique quotient of
the curve $Z$ which identifies the pairs of the
$l_{ij}$s. The curve $C$ is the unique curve whose Jacobian is isomorphic
to $\prym(Z/X)$.

An automorphism of the projective plane $|K_C|$ which fixes the set of points
$l_{ij}$, and their pairs structure, fixes the curve $F$ as
a set. Thus it induces an automorphism of the curve $Z$. this automorphism
commutes with the projection to $Z\ra X$.
Whence it induces a automorphism of
\[
  \prym(Z/X)\cong\jac(C).
\]
This automorphism induces an automorphism of image of $C$ under the canonical
embedding.
\end{proof}
%
\section{Reconstruction of smooth plane quartics form their bitangents}\label{SCS}
%
\begin{dsc}
We consider only non-hyperelliptic curves $C$.
By Theorem \ref{TP2prop}, Proposition \ref{Pauto} and Lemma \ref{Lsingrecon},
we know how to
reconstruct a quartic given the bitangents plus
some combinatorial structure (the 12-tuples, and pairing on each of them).
We have to relax this assumption. The technique we use is group theoretic.
\end{dsc}
\begin{ntt}\label{Ncomb}
Any element $\gc\in\jac(C)[2]\sm\{0\}$ marks a 12-tuple on the 28 bitangents.
Denote this 12-tuple by $\gth(\gc)$.
These 12-tuples has some combinatorial intersection structure, which we call
{\em the} combinatorial structure associated to the 28 bitangents. We call
the pairing induced by $+\gc$ on the set $\gth(\gc)$ {\em the} pairing on $\gth(\gc)$.
\end{ntt} 
\begin{prp}\label{Psymp}
The image of the map
\[
  \begin{aligned}
  \text{pairs of points in }\gS/\ga&\ra\jac(C)[2]/\ga\\
   \{q_i,q_j\}   &\mapsto q_i-q_j
  \end{aligned}
\]
is an isomorphism to $(\ga^\perp/\ga)\sm\{0\}$.
The pullback to of the Weil pairing on the R.H.S. under this isomorphism
is given by the intersection pairing.
\end{prp}
\begin{proof}
Denote by $G_\ga$ the stabilizer of $\ga$ in $SP_6(2)$.
Note that:
\begin{itemize}
\item The map above is a non-trivial $G_\ga$ equivariant map
\item Since $SP_6(2)$ is transitive on
pairs of element in $\jac(C)[2]\sm\{0\}$ with the same
pairing, the action of $G_\ga$ on $\jac(C)[2]/\ga$ has exactly 3 orbits:
\[
  0,(\ga^\perp/\ga)\sm\{0\},(\jac(C)[2]\sm\ga^\perp)/\ga
\]
of sizes $1,15,16$ respectively.
\end{itemize}
Since $6\cdot5/2=15$, the image is the second orbit. The action of the group
$G_\ga$ on both sides map pairs with the same intersection
cardinality to elements with the same Weil pairing. By another counting
argument, the Weil pairing is $1$ if the two sets on the L.H.S. intersect in
a non trivial set.
\end{proof}
\begin{cor}\label{Clev2}
The combinatorial structure on the 28 bitangents determines the pairing on each
of the $\gth(\gc)$.
\end{cor}
\begin{proof}
Denote by $\gc,\gc'$ two distinct elements in $\jac(C)[2]\sm\{0\}$.
\item{\em Claim 1} If $\langle \gc,\gc'\rangle=0$ then
$\gth(\gc)\cap\gth(\gc')$ is a translate of the isotropic group
$\{0,\gc,\gc',\gc+\gc'\}$: 
Note that if $\gth_1\in\gth(\gc)\cap\gth(\gc')$ then
\[
 \gth_1+\{0,\gc,\gc',\gc+\gc'\}\subset\gth(\gc)\cap\gth(\gc').
\]
However, by Proposition \ref{Psymp} there is exactly one translate of the
group $\{0,\gc,\gc',\gc+\gc'\}/\gc$ in $\gth(\gc)/\gc$. Thus, there is
exactly one translate of  $\{0,\gc,\gc',\gc+\gc'\}$ in $\gth(\gc)$.
\item{\em Claim 2}  If $\langle \gc,\gc'\rangle=1$ then
$3|\#\gth(\gc)\cap\gth(\gc')$:
Since the group $SP_6(2)$ acts transitively on pairs of points in
$\jac(C)[2]\sm\{0\}$ with the same Weil pairing the number
$\#\gth(\gc)\cap\gth(\gc')$ is constant. 
By a similar transitivity argument for each $p\in\gth(\gc)$ then number
\[
  k:=\#\{\gc'\in\jac(C)[2]\sm\gc^\perp|p\in\gth(\gc')\}
\]
is constant. Whence,
\[
63\cdot30\cdot \#\gth(\gc)\cap\gth(\gc')=63\cdot 12\cdot k\Rightarrow
  3|\#\gth(\gc)\cap\gth(\gc')
\]
\item{\em Claim 3} The combinatorial structure determines the pairing
on the 12-tuples: for any $\gc\in\jac(C)[2]\sm\{0\}$ the equality
\[
  \gc^\perp\sm\{\gc,0\}=\{\gc'|\#\gth(\gc')\cap\gth(\gc)=4\}.
\]
holds. Moreover the set $\{q_1,\ldots q_6\}$ is the set of non trivial
intersections of the sets $\gth(\gc)\cap\gth(\gc')$ where $\gc'\in\gc^\perp\sm\{\gc,0\}$.
\end{proof}
\begin{prp}
There exists a unique (up to automorphisms of $S_{28}$) embedding
\[
  SP_6(2)\hra S_{28}.
\]
\end{prp}
\begin{proof}
See \cite{DM} chapter 7.7 for proof of the uniqueness. see \cite{DM} p. 46,
\cite{Har} p. 705 or \cite{CCNPW} p. 46 for the construction of the embedding.
\end{proof}
\begin{prp}\label{PsubS28}
The only groups that strictly contain $SP_6(2)$ under the embedding
$SP_6(2)\hra S_{28}$, are $A_{28},S_{28}$.
\end{prp}
\begin{proof}
There are only 6 2-transitive groups on 28 elements (see \cite{DM} chapter 7.7
for a complete list of the finite 2-transitive groups)
\[
  \begin{array}{ll}
  \mathrm{Group} & \mathrm{order} \\
\hline
  S_{28}     & 28!             \\
  A_{28}     & 28!/2           \\ 
  SP_6(2) & 63\cdot 30\cdot 12 \cdot 8^2=1451520 \\
  PSL_2(9) & \frac{(9^2-1)(9^2-9)}{2(9-1)}=360  \\
  PSU_3(3) & 6048  \\ 
  R(3)     & (3^3+1)3^3(3-1)=1512.
  \end{array}
\]
(See \cite{CCNPW} p. x for the size of $PSU_3(3)$ and \cite{DM} p. 251 for the
formula of the size of $R(q)$.)
\end{proof}
\begin{proof}[proof of Theorem \ref{Tbitan}]
We consider the following objects:
\begin{description}
\item[$\fC$] The set of smooth plane quartics.
\item[$\ti{\fB}$] The set of possible sets of bitangents + the combinatorial
   structure of the 12-tuples.
\item[$\fB\subset\Sym^{28}\BP^2$] The set of possible sets of bitangents.
\item[$\ol{\fB}\subset(\BP^2)^{28}$] The pullback of $\fB$ under
     $(\BP^2)^{28}\ra\Sym^{28}\BP^2$.
\item[$\ol{\fB}_0\subset\ol{\fB}$] a connected component (they are
   all isomorphic).
\end{description}
Note that $C$ is an irreducible quasi-projective variety.
The cover $\ol{\fB}/\fB$ is etale and Galois, therefore so is
$\ol{\fB}_0/\fB$.
By Corollary \ref{Clev2}, the map $\fC\ra\ti{\fB}$ is 1-1.
By \cite{Har} Result (ii) p. 687, the monodromy group of the 28 bitangents:
$Gal(\ol{\fB}_0/\fC)$, is $SP_6(2)$.
By definition,
\[
  Gal(\ol{\fB}_0/\ti{\fB})\subset Gal((\BP^2)^{28}/\Sym^{28}\BP^2)=S_{28}.
\]
The situation can be summarized in the following diagram:
\[
  \xymatrix{
  & \ol{\fB}_0 \ar@{^{(}->}[r]\ar[d] & (\BP^2)^{28} \ar[dd] \\
  & \ti{\fB} \ar[d] \ar@/^/@{.}[u]^{SP_6(2)}\\
  \fC \ar[ur]^{\cong} \ar[r] & \fB \ar@{^{(}->}[r] \ar@{.}@/_/[uu]_{\text{?}} &
  \Sym^{28}\BP^2 \ar@{.}@/_/[uu]_{S_{28}}
}
\] 
Let us suppose that $Gal(\ol{\fB}_0/\fB)\varsupsetneq SP_6(2)$. By Proposition
\ref{PsubS28}, $Gal(\ol{\fB}_0/\fB)\supset A_{28}$.
Consider a 12-tuple of the associated combinatorial level-2 structure.
Denote by $F$ the unique cubic on which these 12 bitangents sit.
By Lemma \ref{LWirr}, if the curve $F$ is not smooth, it is reducible.
In this case there
is no 3-cycle which permutes the $l_{ij}$s between the irreducible components
of $F$, but fixes all but one of the pairs in each component.
If the curve $F$ is smooth then
by Lemma \ref{Clev2} and Proposition \ref{P9pts} there is no
3-cycle in $Gal(\ol{\fB}_0/\fB)$ that breaks the unique level pairing on
the 12-tuple.
\end{proof}
\begin{rmr}\label{RaCS}
It seems that an alternative proof for the reconstructability can be given
using Aronhold systems, along the same lines of the proof above:
\item Aronhold systems are preserved by the Galois group action on the
bitangents. Let us find the group $\ti{G}\subset SP_6(2)$ fixing a specific system.
Since $SP_6(2)$ is transitive on bitangents,
$\ti{G}\subsetneq SP_6(2)$. Whence, if $K\subset SP_6(2)$ is the
maximal subgroup containing $\ti{G}$, then $[K:\ti{G}]|288$. There is only one such
group: $O_6^+$, where $[SP_6(2):O_6^+]=36$, and $O_6^+\cong S_8$. This group
is also transitive on the bitangents, so $\ti{G}\subsetneq O_6^+$. The group
$S_8$ has only two subgroups $\ol{G}$ s.t. $[S_8:\ol{G}]|8$: $A_8$ and $S_7$. The
group $A_8$ is again ``too transitive'' (since $[S_8:A_8]=2$), 
so $\ti{G}=S_7$.
The kernel of the action of $S_7$ on the bitangents is either one
of $1,A_7,S_7$.
As $O_6^+$ acts transitively on the bitangents, and $[O_6^+:\ti{G}]=8$, the kernel
has to be $1$. i.e.:
the action of $\ti{G}$ on the system is the symmetric action.
We use a similar technique to the one applied in the proof of theorem
\ref{Tbitan}:
The subgroup of $A_{28}$ acting on a 7-tuple via $S_7$ is
$A_{28}\cap(S_{21}\times S_7)$. The group $\ti{G}$ is a strict subgroup of this group.
\item The details of this sketch have not been thoroughly checked.
\end{rmr}


\end{document}